\documentclass[12pt]{article}
\usepackage{graphicx}
\usepackage{geometry}
\usepackage{relsize}
\author{\relsize{+1}John Lindgren and Vibeke Libby\\ \relsize{-1}\sl john.lindgren@stanfordalumni.org\\  \relsize{-1}\normalfont Agilemath}
\title{\bf An Arithmetic and Geometric Mean Invariant}
\begin{document}
\maketitle
\begin{abstract}
A positive real interval, $[a, b]$ can be partitioned into sub-intervals such that sub-interval widths divided by sub-interval "`average"' values remains constant. That both Arithmetic Mean and Geometric Mean "`average"' values produce constant ratios for the same log scale is the stated invariance proved in this short note. The continuous analog is briefly considered and shown to have similar properties.
\end{abstract}

\section{Problem and Summary}

Define a relative weight for an interval, $[a, b]$, with $0 < a < b$, as the interval length divided by the average value in the interval. We explore interval partitions of [a, b] that have equal relative weights for consecutive sub-intervals and show that there is a surprising invariance to whether we use arithmetic or geometric averages, or even an interval endpoint, for the definition of "`average."' The ratios shown in figure 1 for a partition of $[1, 16]$ are seen to be constant when using either the arithmetic or geometric means of sub-interval endpoints. 

\section{Main Result}
\bf {Theorem 1:} Given an interval $[a, b]$, there is a partition
\begin{equation}
a = x_0 < x_1 < ... < x_n = b
\end{equation}
such that the following equations are satisfied:

\begin{equation}
{\frac{x_{i-1} - x_i}{(x_{i+1} + x_i)/2}} = {\frac{x_{j+1} - x_j}{(x_{j+1} + x_j)/2}}; \forall i < j \in \{0,1...n\}
\end{equation}
\\

\begin{equation}
{\frac{x_{i+1} - x_i}{\sqrt{x_{i+1}x_i}}} = {\frac{x_{j+1} - x_j}{\sqrt{x_{j+1}x_j}}}; \forall i < j \in \{0,1...n\}
\end{equation}
\\

Moreover, the partition satisfying the conditions is given by:\\

\begin{equation}
x_i = b^{i/n}a^{(n-i)/n}
\end{equation}

\begin{figure}[ht!]
\vspace*{-20.ex}
    \begin{center}
            \includegraphics[width=1.1\textwidth,height=1.1\textheight]{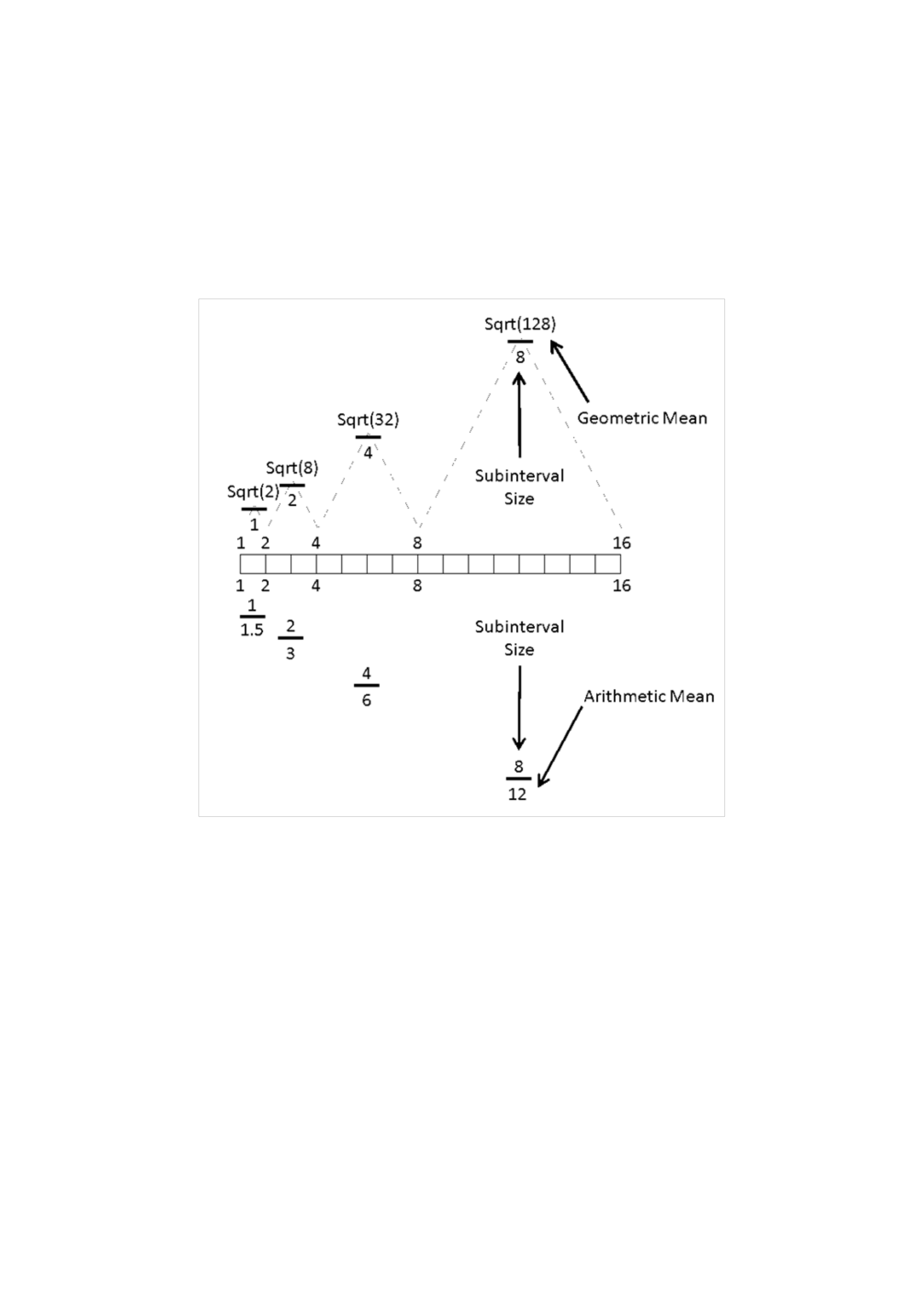}\\
\vspace*{-40.ex}
    \caption{\it Interval endpoints plotted on a log scale for [1, 16]  with n=4 sub-intervals.}
    \end{center}
\vspace*{-2.ex}
\end{figure}

\section{Proof}
The invariance follows by considering the ratio between consecutive terms:

\begin{equation}
z_{i+1} = \frac{x_{i+1}}{x_i}
\end{equation}

\subsection{The AM case}
The arithmetic mean (AM) case follows from:

\begin{equation}
{\frac{z_{i+1}-1}{z_{i+1}+1}} = {\frac{z_i-1}{z_i+1}}
\end{equation}
Reduction yields:
\begin{equation}
(z_{i+1}-1)(z_i + 1) = (z_{i+1} + 1)(z_i - 1) \Rightarrow
\end{equation}

\begin{equation}
{z_{i+1} - z_i - 1 = - z_{i+1} + z_i -1} \Rightarrow z_{i+1} = z_i
\end{equation}

\subsection{The GM case}
The geometric mean (GM) case follows similarly:
\begin{equation}
{\frac{z_{i+1}-1}{\sqrt{x_{i+1}x_i}/x_i}} = {\frac{z_i-1}{\sqrt{x_ix_{i-1}}/x_{i-1}}} = {\frac{z_{i+1}-1}{\sqrt{z_{i+1}}}} = \frac{z_i - 1}{\sqrt{z_i}} \Rightarrow
\end{equation}

\begin{equation}
{\sqrt{z_i}(z_{i+1} - 1)} = {\sqrt{z_{i+1}}(z_i - 1)} \Rightarrow 
\end{equation}
\\
\begin{equation}
\sqrt{z_iz_{i+1}}\Big{(}{\sqrt{z_{i+1}}} - 1/{\sqrt{z_iz_{i+1}}}\Big{)} - {\sqrt{z_iz_{i+1}}}\Big{(}{\sqrt{z_i}} - 1/{\sqrt{z_iz_{i+1}}}\Big{)}=0 
\end{equation}
\\
Factoring the last equation results in:
\begin{equation}
\sqrt{z_iz_{i+1}}\Big{(}\sqrt{z_{i+1}} - \sqrt{z_i}\Big{)} = 0 \Rightarrow z_{i+1} = z_i
\end{equation}

Finally, the actual partition can now be determined for both the AM and the GM cases:
\begin{equation}
z_{i+1} = \frac{x_{i+1}}{x_i} = c
\end{equation}

\begin{equation}
b = x_n = cx_{n-1} = ... = c^nx_0 = c^na
\end{equation}

\begin{equation}
c = (b/a)^{1/n} \Rightarrow
\end{equation} 

\begin{equation}
x_i = cx_{i-1} = (b/a)^{1/n}x_{i-1} = (b/a)^{2/n}x_{i-2} = ... = (b/a)^{i/n}x_0 = b^{i/n}a^{(-i)/n}a = b^{i/n}a^{(n-i)/n}
\end{equation}
QED

\section{Continuous Analog Mean}
Using an interval endpoint for the "`average"' of each sub-interval yields exactly the same recurrence relationship. As $n$ goes to infinity, the corresponding differential equation is given by:
\begin{equation}
\frac{ \Delta x}{x} = k\Delta n \Rightarrow
\end{equation}

\begin{equation}
\int\frac { dx}{x} = k\int dn \Rightarrow
\end{equation}
 
\begin{equation}
ln(x) = kn + c \Rightarrow x(n) = Ae^{Bn}
\end{equation}

When the same boundary conditions specified in the discrete case are given for $x$, the constants $A$ and $B$ in the above equation produce the same log scale solution as the discrete cases.

\section{Background}
The problem has its origin in generating test sets for an optical bench. Part of the EM spectrum was to be sampled into bands and the bands were to be selected such that the relative error would remain roughly constant for all bands:

\begin{equation}\frac{\Delta\lambda}{\lambda_{avg}}
\end{equation}

Further analysis showed that the choice of "`average"' made no difference in selecting test bands.

\section{References}
\begin{enumerate}
\item Stanislas Dehaene, Véronique Izard, Elizabeth Spelke, and Pierre Pica. Log or Linear? Distinct Intuitions of the Number Scale in Western and Amazonian Indigene Cultures. Science, 2008; 320 (5880): 1217 DOI: 10.1126/science.1156540

\end{enumerate}
\end{document}